\documentclass{article}
\usepackage{verbatim}
\usepackage{amsmath}
\usepackage[mathscr]{eucal}
\usepackage{amsthm}
\usepackage{amssymb,mathrsfs}
\usepackage[german,english]{babel}
\theoremstyle{plain}
\newtheorem{Theo}{Theorem}[section]
\newtheorem{Prop}[Theo]{Proposition}
\newtheorem{Lem}[Theo]{Lemma}

\newtheorem{IntroTh}{Theorem}
\theoremstyle{definition}

\theoremstyle{remark}

\newcommand{\gl}{\mathfrak{g}\mathfrak{l}}

\makeatletter
	
	\@addtoreset{equation}{section}
\makeatother
\begin{document}
\title{A super Frobenius formula for the characters of Iwahori-Hecke algebras}
\author{Hideo Mitsuhashi\footnote{E-mail\,:\,mitsu@gt.tomakomai-ct.ac.jp}}

\date{}

\maketitle

\begin{center}
Department of Mechanical Engineering \\
Tomakomai National College of Technology \\
443 Nishikioka, Tomakomai, Hokkaido, 059-1275, Japan
\end{center}

\begin{abstract}
In this paper, we establish a super Frobenius formula for the characters of Iwahori-Hecke algebras. 
We define Hall-Littlewood supersymmetric functions in a standard manner to make 
supersymmetric functions from symmetric functions, and investigate some properties of supersymmetric functions. 
Based on Schur-Weyl reciprocity between Iwahori-Hecke algebras and the general quantum 
super algebras, which was obtained in \cite{Mit2}, we derive 
that one of several types of Hall-Littlewood supersymmetric functions, up to constant, 
generates the values of the irreducible characters of Iwahori-Hecke algebras 
at the elements corresponding to cycle permutations. Our formula in this article includes both the ordinary 
quantum case ($n=0$) that was obtained in \cite{Ram} and the classical super case ($q{\rightarrow}1$). 
\end{abstract}

\pagestyle{plain}
\pagenumbering{arabic}
\renewcommand{\theenumi}{\arabic{enumi}}
\renewcommand{\labelenumi}{(\theenumi)}

\section{Introduction}
The Frobenius formula is one of powerful methods to compute the irreducible characters of symmetric 
groups. It is based on Schur-Weyl reciprocity; there are actions of symmetric groups and 
of general linear groups which generate the full centralizer of each other. Schur-Weyl reciprocity has 
been extended to various groups and algebras up to the present. 
Among them, two remarkable extensions for us are the super type extension (\cite{B-R},\cite{Ser}) and 
the quantum type one (\cite{Jimbo}). 
In our paper \cite{Mit2}, we established Schur-Weyl reciprocity between the 
Iwahori-Hecke algebra $\mathscr{H}_{\mathbb{Q}(q),r}(q)$ and the quantum super algebra 
$U_q^\sigma\big{(}{\gl}(m,n)\big{)}$, 
which unifies Schur-Weyl reciprocity of super type and that of quantum type. 
In \cite{Ram}, Ram gave a (ordinary) Frobenius formula for the characters of the Iwahori-Hecke algebras of 
type $A$, which is based on the Schur-Weyl reciprocity between the Iwahori-Hecke algebra of type $A$ and the 
quantum enveloping algebra of $\gl_n$ that was given in \cite{Jimbo}. 
An extension of Frobenius formula to Ariki-Koike algebras, that are Hecke algebras associated to 
complex reflection groups $G(r,1,n)$, is found in \cite{Shoji}. 
In this paper, we give a super Frobenius formula for the characters of the Iwahori-Hecke algebras of 
type $A$ that extends the result of Ram. \par
As in the Macdonald's book \cite{Mac}, symmetric functions play a crucial role in the representation theory 
of symmetric groups. 
Especially, the transition matrix $M(p,s)$ from power sum functions to Schur functions is the character table 
of the symmetric group. Combinatorial rules to compute character values such as Murnaghan-Nakayama formula 
have been given by making use of properties of symmetric functions. 
On the other hand, various generalizations of symmetric functions are defined until now. 
Hall-Littlewood functions are symmetric functions with one parameter which intermediate between 
monomial symmetric functions and Schur functions. 
Among several types of Hall-Littlewood functions, The one which is denoted by $q_{\lambda}(x,q)$ 
yields the Frobenius formula for the characters of the Iwahori-Hecke algebra (\cite{Ram}). \par
In order to give a super Frobenius formula for the Iwahori-Hecke algebra, we define the Hall-Littlewood 
supersymmetric functions $P_{\lambda}(x/y;q)$ and $q_{\lambda}(x/y,q)$ and investigate some properties 
of them in section 2. 
Let $x=(x_1,x_2,{\hdots},x_m)$ and $y=(y_1,y_2,{\hdots},y_n)$ be commutative variables. 
Using a partition of unity in $\mathscr{H}_{\mathbb{Q}(q),r}(q)$, that is a complete set 
of orthogonal minimal idempotents which is specialized to a partition of unity in $\mathbb{C}[\mathfrak{S}_r]$, 
we derive that the trace of the product of $\pi_r(h){\in}\mathcal{A}_q$ 
and $D_r{\in}\mathcal{B}_q$ is R.H.S. of the super Frobenius formula (Definitions of $\pi_r,\mathcal{A}_q,
D_r,\mathcal{B}_q$ are given in section 3). 
\renewcommand{\theIntroTh}{4.4}
\begin{IntroTh}
For any $h{\in}\mathscr{H}_{\mathbb{Q}(q)(x,y),r}(q)$, 
\[
\operatorname{tr}(D_r{\pi}_r(h))=\sum_{\lambda{\vdash}r}\chi^{\lambda}(h)s_{\lambda}(x/y),
\]
where $\chi^{\lambda}$ is the irreducible character of $\mathscr{H}_{\mathbb{Q}(q)(x,y),r}(q)$ 
corresponding to $\lambda$. 
\end{IntroTh}
We investigate L.H.S. of Theorem 4.4 in section 5. We obtain that when $h$ is the element $T_{\gamma_k}$ 
of $\mathscr{H}_{\mathbb{Q}(q)(x,y),r}(q)$ corresponding to 
the cycle permutation of length $k$, L.H.S. coincides with $q_k(x/y;q^{-2})$ up to constant. 
\renewcommand{\theIntroTh}{5.3}
\begin{IntroTh}
\[
\operatorname{tr}\big{(}D_k{\pi_k}(T_{\gamma_k})\big{)}=\dfrac{q^k}{q-q^{-1}}q_k(x/y;q^{-2}).
\]
\end{IntroTh}
Finally, for the product $T_{\gamma_{\mu}}=T_{\gamma_{\mu_1}}T_{\gamma_{\mu_2}}{\cdots}
T_{\gamma_{\mu_l}}$ where $\mu=(\mu_1,\mu_2,{\hdots}){\vdash}r$, we have 
\renewcommand{\theIntroTh}{5.5}
\begin{IntroTh}
For $\mu{\vdash}r$, 
\[
\dfrac{q^{|{\mu}|}}{(q-q^{-1})^{l({\mu})}}q_{\mu}(x/y;q^{-2})
=\sum_{\lambda{\vdash}r}\chi^{\lambda}(T_{\gamma_{\mu}})s_{\lambda}(x/y).
\]
\end{IntroTh}
Though not all the values of characters can be computed by Theorem 5.5, owing to our version of 
Ram's result 
\renewcommand{\theIntroTh}{5.6}
\begin{IntroTh}[\cite{Ram}, Theorem 5.1]
For each $T_{\sigma},\,\sigma{\in}\mathfrak{S}_r$, there exists a $\mathbb{Z}[q,q^{-1}]$ linear combination 
\[
c_{\sigma}=\sum_{\mu{\vdash}r}a_{{\sigma}{\mu}}T_{\gamma_{\mu}},
\]
$a_{{\sigma}{\mu}}{\in}\mathbb{Z}[q,q^{-1}]$, such that $\chi(T_{\sigma})=\chi(c_{\sigma})$ 
for all characters $\chi$ of $\mathscr{H}_{\mathbb{Q}(q),r}(q)$,
\end{IntroTh}
shows that any character of $\mathscr{H}_{\mathbb{Q}(q),r}(q)$ is determined by its values 
on the elements $T_{\gamma_{\mu}}$. \par
Our result extends \cite{Ram}, and may yield new combinatorial rules for computing 
the irreducible characters of Iwahori-Hecke algebras. 
Moreover, as we have pointed out in \cite{Mit2}, the super version of the representation theory of 
the symmetric group and the Iwahori-Hecke algebra are more suitable to describe the representation theory of 
the alternating group and its $q$-analogue. 
Our result will be used to derive a Frobenius formula for the $q$-analogue of the alternating group which 
was defined in \cite{Mit}. 
\section{Preliminaries on supersymmetric functions}
Symmetric functions, besides their own interest, play important roles in various areas in mathematics. 
Particularly, the relation between symmetric functions and the representation theory is intimate. 
In this section, we will give supersymmetrizations of various classes of symmetric functions and 
investigate relationships between them. The basic reference of symmetric functions and supersymmetric 
functions is Macdonald's book \cite{Mac}. We will follow \cite{Mac} with respect to our notation about 
symmetric functions unless otherwise stated. \par
We denote by $\mathfrak{S}_m$ the symmetric group of order $m!$. 
Let $\Lambda_m=\mathbb{Z}[x_1,x_2,{\hdots},x_m]^{\mathfrak{S}_m}$ be the ring of symmetric functions 
of $m$ variables and $(\Lambda_m)_{\mathbb{Q}}=\Lambda_m{\otimes_{\mathbb{Z}}\mathbb{Q}}$.
One can define an involution 
\[
\omega:\Lambda_m{\longrightarrow}\Lambda_m
\]
by $\omega(e_r)=h_r$. Then, as \cite{Mac},I,(2.13),(3.8), 
$\omega(p_r)=(-1)^{r-1}p_r$ and $\omega(s_{\lambda})=s_{\lambda'}$ hold. 
Let $x=(x_1,x_2,{\hdots}x_m)$ be $m$ variables and $y=(y_1,y_2,{\hdots}y_n)$ $n$ variables. 
We denote by $f(x,y)$ a symmetric function $f(x_1,x_2,{\hdots},x_m,y_1,y_2,{\hdots},y_n)
{\in}\Lambda_{m+n}$ for abbreviation. 
According to \cite{Mac},I.3,ex.$23$, we set 
\begin{equation}
E_{x/y}(t)=\prod_{i=1}^m(1+x_it)\Big{/}\prod_{j=1}^n(1+y_jt).
\end{equation}
The elementary supersymmetric functions $e_k(x/y)\,(k{\geq}0)$ is defined by 
\begin{equation}
E_{x/y}(t)=\sum_{k{\geq}0}e_k(x/y)t^k.
\end{equation}
From (2.1),(2.2), $e_k(x/y)$ is a polynomial in $x=(x_1,x_2,{\hdots}x_m)$ and $y=(y_1,y_2,{\hdots},y_n)$, 
symmetric in each set of variables separately. 
We denote by $\Lambda_{m,n}$ the ring of polynomials in $x_1,x_2,{\hdots}x_m,y_1,y_2,{\hdots}y_n$, 
which are separately symmetric in $x's$ and $y's$, namely 
\[
\Lambda_{m,n}=\mathbb{Z}[x_1,x_2,{\hdots},x_m]^{\mathfrak{S}_m}
{\otimes_{\mathbb{Z}}}\mathbb{Z}[y_1,y_2,{\hdots},y_n]^{\mathfrak{S}_n}.
\]
We also denote by $(\Lambda_{m,n})_{\mathbb{Q}}=\Lambda_{m,n}{\otimes_{\mathbb{Z}}\mathbb{Q}}$. 
$\omega$ can be regarded as an involution on $\Lambda_{m,n}$ which acts on each factor of 
$\mathbb{Z}[x_1,x_2,{\hdots},x_m]^{\mathfrak{S}_m}
{\otimes_{\mathbb{Z}}}\mathbb{Z}[y_1,y_2,{\hdots},y_n]^{\mathfrak{S}_n}$ as $\omega$. 
One can deduce from (2.1),(2.2) that 
\begin{equation*}
\begin{split}
e_k(x/y)&=\sum_{i+j=k}(-1)^je_i(x)h_j(y){\quad}\text{(\cite{Mac},$\rm{I}.3, ex.23$)}\\
&=\sum_{i+j=k}e_i(x)h_j(-y)
\end{split}
\end{equation*}
Let $\omega_y$ be the involution $\omega$ acting on symmetric functions of $y$ variables. 
Then, $e_k(x/y)$ is obtained from ${\omega_y}e_k(x,y)$ by replacing each $y_i$ by $-y_i$. 
For each partition $\lambda$, supersymmetric Schur function $s_{\lambda}(x/y)$ is defined by 
\begin{equation*}
s_{\lambda}(x/y)={\det}\big{(}e_{\lambda'_i-i+j}(x/y)\big{)}_{1{\leq}i,j{\leq}l}
\end{equation*}
where $l=l(\lambda')$. 
We readily see that 
\begin{equation*}
e_k(x/y)=s_{(1^k)}(x/y)
\end{equation*}
By making use of the equations 
\begin{equation*}
{\omega}s_{\lambda/\mu}(x)=s_{\lambda'/\mu'}(x)
{\quad}\text{(\cite{Mac},I,$(5.6)$)},
\end{equation*}
\begin{equation*}
s_{\lambda}(x,y)=\sum_{\mu{\subset}\lambda}s_{\mu}(x)s_{{\lambda}/{\mu}}(y) 
{\quad}\text{(\cite{Mac},I,$(5.9)$)},
\end{equation*}
and $\omega_y$, one can obtain 
\begin{equation}
\begin{split}
s_{\lambda}(x/y)&=\sum_{\mu{\subset}\lambda}(-1)^{|\lambda-\mu|}s_{\mu}(x)s_{{\lambda'}/{\mu'}}(y) 
{\quad}\text{(\cite{Mac},$\rm{I}.5, ex.23(1)$)}\\
&=\sum_{\mu{\subset}\lambda}s_{\mu}(x)s_{{\lambda'}/{\mu'}}(-y).
\end{split}
\end{equation}
In a similar manner to (2.1), we define $H_{x/y}(t)$ to be 
\begin{equation}
H_{x/y}(t)=\prod_{i=1}^m(1-x_it)^{-1}\prod_{j=1}^n(1-y_jt)
\end{equation}
and the completely supersymmetric functions $h_k(x/y)\,(k{\geq}0)$ to be 
\begin{equation}
H_{x/y}(t)=\sum_{k{\geq}0}h_k(x/y)t^k.
\end{equation}
One can deduce from (2.4),(2.5) that 
\begin{equation*}
h_k(x/y)=\sum_{i+j=k}h_i(x)e_j(-y)=e_k(-y/-x).
\end{equation*}
Moreover, applying $\omega$ to $h_k(x/y)$, we have
\begin{equation*}
{\omega}h_k(x/y)=\sum_{i+j=k}{\omega}h_i(x){\omega}e_j(-y)=\sum_{i+j=k}e_i(x)h_j(-y)=e_k(x/y)
=s_{(1^k)}(x/y).
\end{equation*}
Hence we obtain 
\begin{equation*}
h_k(x/y)={\omega}s_{(1^k)}(x/y)=s_{(k)}(x/y).
\end{equation*}
We may also define the power sum supersymmetric functions 
$p_k(x/y)\,(k{\geq}1)$ as follows. 
\begin{equation*}
\begin{split}
p_k(x/y)&={\omega_y}p_k(x,-y)\\
&=p_k(x)+(-1)^{k-1}p_k(-y)\\
&=p_k(x)-p_k(y)
\end{split}
\end{equation*}
In the same way as \cite{Mac},I,(2.10), we obtain the following relation 
for the generating function $P_{x/y}(t)$. 
\begin{equation*}
\begin{split}
P_{x/y}(t)&=\sum_{r{\geq}1}p_r(x/y)t^{r-1}\\
&=\sum_{r{\geq}1}\sum_{i=1}^mx_i^rt^{r-1}-\sum_{r{\geq}1}\sum_{i=1}^ny_i^rt^{r-1}\\
&=\sum_{i=1}^m\dfrac{x_i}{1-x_it}-\sum_{i=1}^n\dfrac{y_i}{1-y_it}\\
&=\sum_{i=1}^m\dfrac{d}{dt}{\log}\dfrac{1}{1-x_it}-\sum_{i=1}^n\dfrac{d}{dt}{\log}\dfrac{1}{1-y_it}\\
&=\dfrac{d}{dt}{\log}\prod_{i=1}^m(1-x_it)^{-1}\prod_{j=1}^n(1-y_jt)\\
&=\dfrac{d}{dt}{\log}H_{x/y}(t)=\dfrac{H_{x/y}'(t)}{H_{x/y}(t)}.
\end{split}
\end{equation*}
Thus, we have 
\begin{equation*}
\begin{split}
H_{x/y}(t)&={\exp}\big{(}\sum_{r{\geq}1}p_r(x/y)\dfrac{t^r}{r}\big{)} \\
&=\prod_{r{\geq}1}{\exp}\big{(}p_r(x/y)\dfrac{t^r}{r}\big{)} \\
&=\prod_{r{\geq}1}\sum_{m_r{\geq}0}\dfrac{\big{(}p_r(x/y)t^r\big{)}^{m_r}}{m_r!r^{m_r}},
\end{split}
\end{equation*}
where $m_i$ is the number of appearances of $i$ in 
$\lambda_1,\lambda_2,{\hdots},\lambda_{l(\lambda)}$. 
Let $\lambda?=\prod_{i{\geq}1}i^{m_i}m_i!$. Then 
\begin{equation*}
\begin{split}
H_{x/y}(t)=\sum_{\lambda}\dfrac{p_{\lambda}(x/y)}{\lambda?}t^{|\lambda|} 
\end{split}
\end{equation*}
summed over all partitions $\lambda$. On the other hand, 
\[
H_{x/y}(t)=\sum_{k{\geq}0}h_k(x/y)t^k,
\]
so we can conclude that 
\begin{equation}
h_k(x/y)=\sum_{\lambda{\vdash}k}\dfrac{p_{\lambda}(x/y)}{\lambda?}.
\end{equation}
We notice $e_r(x/y)$ and $h_r(x/y)$ are defined to be zero if $r<0$. \par
By an analogy of the case of symmetric functions, 
we will use the following equation later to give a super version of the Frobenius formula 
for the characters of Iwahori-Hecke algebras. 
\begin{Theo}[Super Jacobi-Trudi formula,\cite{Remmel},(1.9)]
Let $\lambda{\vdash}r$. Then 
\begin{equation*}
s_{\lambda}(x/y)=\det[h_{\lambda_i-i+j}(x/y)]_{1{\leq}i,j{\leq}r}
\end{equation*}
\end{Theo}
In \cite{Remmel}, the above formula is given for the hook-Schur function 
\[
HS_{\lambda}(x;y)=\sum_{\mu{\subset}\lambda}s_{\mu}(x)s_{\lambda'/\mu'}(y),
\]
which is different from $s_{\lambda}(x/y)$ slightly. But replacing $y$ with $-y$, Theorem 2.1 is 
obtained from \cite{Remmel} immediately. 
Super Jacobi-Trudi formula is also found in \cite{P-H}. \par
Let $\Lambda_{m,n}^r$ be the subspace of $\Lambda_{m,n}$ consisting of $r$ homogeneous polynomials. 
For arbitrary $r{\in}\mathbb{N}$, we define a map $\psi$ from $\mathfrak{S}_r$ to $\Lambda_{m,n}^r$ 
as follows. 
\begin{equation*}
\psi(w)=p_{\rho(w)}(x/y)=p_{\rho_1}(x/y)p_{\rho_2}(x/y){\cdots}p_{\rho_{l(\rho)}}(x/y)
\end{equation*}
where $\rho(w)=(\rho_1,\rho_2,{\hdots}\rho_{l(\rho)})$ is the cycle type of $w$. 
Let $v{\times}w{\in}\mathfrak{S}_{r_1}{\times}\mathfrak{S}_{r_2}$. 
Then $v{\times}w$ can be embedded in $\mathfrak{S}_{r_1+r_2}$ in many ways. 
Nevertheless, for the same reason as in $\rm{I}.7$ in \cite{Mac}, we have 
\begin{equation*}
\psi(v{\times}w)=\psi(v)\psi(w).
\end{equation*}
In general, for functions $f,g$ on a finite group $G$ which values in a commutative $\mathbb{Q}$-algebra, 
the scalar product of $f$ and $g$ is defined by 
\[
{\langle}f,g{\rangle}=\dfrac{1}{|G|}\sum_{x{\in}G}f(x)g(x^{-1}).
\]
Let $\operatorname{Irr}(G)$ be a basic set of irreducible $\mathbb{C}$-characters of 
$G$ and 
$\operatorname{ch}(\mathbb{C}[G])={\oplus}_{\zeta{\in}\operatorname{Irr}}\mathbb{Z}\zeta$ the 
$\mathbb{Z}$-module of 
virtual characters. 
For $f_i{\in}\operatorname{Irr}(\mathfrak{S}_{r_i})\,(i=1,2)$, we define 
$f_1{\cdot}f_2{\in}\operatorname{ch}(\mathbb{C}[\mathfrak{S}_{r_1+r_2}])$ to be  
\begin{equation}
f_1{\cdot}f_2=\operatorname{ind}_{\mathfrak{S}_{r_1}{\times}\mathfrak{S}_{r_2}}^{\mathfrak{S}_{r_1+r_2}}
(f_1{\times}f_2)
\end{equation}
Then, as in \cite{Mac},I.7, the $\mathbb{Z}$-module 
\[
\operatorname{ch}(\mathbb{C}[\mathfrak{S}])=\bigoplus_{r{\geq}0}\operatorname{ch}(\mathbb{C}[\mathfrak{S}_r])
\]
has a ring structure with the product (2.7). 
By an analogy of the characteristic map as in \cite{Mac},I.7, we define a $\mathbb{Z}$-linear map  
\[
\operatorname{sch} : \operatorname{ch}(\mathbb{C}[\mathfrak{S}]){\longrightarrow}
(\Lambda_{m,n})_{\mathbb{C}}=\Lambda_{m,n}{\otimes}\mathbb{C}
\]
to be 
\begin{equation*}
\operatorname{sch}(f)={\langle}f,\psi{\rangle}_{\mathfrak{S}_r}
=\dfrac{1}{r!}\sum_{w{\in}\mathfrak{S}_r}f(w){\psi}(w)
=\sum_{\mu{\vdash}r}f(\mu)\dfrac{p_{\mu}(x/y)}{{\mu}?},
\end{equation*}
where $f(\mu)$ is the value of $f$ at elements of cycle type $\mu$. 
Using Frobenius reciprocity, we obtain 
\begin{equation*}
\begin{split}
\operatorname{sch}(f_1{\cdot}f_2)&={\langle}\operatorname{ind}_{\mathfrak{S}_{r_1}{\times}\mathfrak{S}_{r_2}}^{\mathfrak{S}_{r_1+r_2}}(f_1{\times}f_2),\psi{\rangle}_{\mathfrak{S}_{r_1+r_2}} \\
&={\langle}(f_1{\times}f_2),
\operatorname{res}_{\mathfrak{S}_{r_1}{\times}\mathfrak{S}_{r_2}}^{\mathfrak{S}_{r_1+r_2}}(\psi){\rangle}
_{\mathfrak{S}_{r_1}{\times}\mathfrak{S}_{r_2}} \\
&=\dfrac{1}{|\mathfrak{S}_{r_1}||\mathfrak{S}_{r_2}|}\sum_{v{\times}w{\in}\mathfrak{S}_{r_1}
{\times}\mathfrak{S}_{r_2}}(f_1{\times}f_2)(v{\times}w)\psi(v{\times}w) \\
&=\dfrac{1}{|\mathfrak{S}_{r_1}||\mathfrak{S}_{r_2}|}\sum_{v{\in}\mathfrak{S}_{r_1},w{\in}\mathfrak{S}_{r_2}}
f_1(v)f_2(w)\psi(v)\psi(w) \\
&={\langle}f_1,\psi{\rangle}_{\mathfrak{S}_{r_1}}{\langle}f_2,\psi{\rangle}_{\mathfrak{S}_{r_2}} \\
&=\operatorname{sch}(f_1)\operatorname{sch}(f_2).
\end{split}
\end{equation*}
Therefore $\operatorname{sch}$ is a ring homomorphism. 
\begin{Prop}
\[
s_{\lambda}(x/y)=\sum_{\mu{\vdash}r}{\chi^{\lambda}(\mu)}\dfrac{p_{\mu}(x/y)}{{\mu}?},
\]
where $\chi^{\lambda}$ is the irreducible character of $\mathbb{C}[\mathfrak{S}_r]$ corresponding to $\lambda$. 
\end{Prop}
\begin{proof}
Let $\eta_r$ be the identity character of $\mathfrak{S}_r$ and $\eta_r=0$ if $r<0$. 
Then by (2.6),
\[
\operatorname{sch}(\eta_r)=\sum_{\mu{\vdash}r}\dfrac{p_{\mu}(x/y)}{\mu?}
=h_r(x/y).
\]
For $\lambda=(\lambda_1,\lambda_2,{\hdots},\lambda_l){\vdash}r\,(l=l(\lambda))$, 
let $\eta_{\lambda}=\eta_{\lambda_1}{\cdot}\eta_{\lambda_2}{\cdot}
{\cdots}{\cdot}\eta_{\lambda_l}$ be the character 
\[
\eta_{\lambda}=\operatorname{ind}_{\mathfrak{S}_{\lambda_1}{\times}
\mathfrak{S}_{\lambda_1}{\times}{\cdots}{\times}\mathfrak{S}_{\lambda_l}}^{\mathfrak{S}_r}
(\eta_{\lambda_1}{\times}\eta_{\lambda_2}{\times}{\cdots}{\times}\eta_{\lambda_{\lambda_l}})
\]
of $\mathfrak{S}_r$. 
Then we have $\operatorname{sch}(\eta_{\lambda})=h_{\lambda}(x/y)$. 
From \cite{Mac},I,(7.4),(7.6)(i), $\chi^{\lambda}={\det}[\eta_{\lambda_i-i+j}]_{1{\leq}i,j{\leq}r}$. 
Using Theorem 2.1 one can deduce 
\[
\operatorname{sch}(\chi^{\lambda})={\det}[h_{\lambda_i-i+j}(x/y)]_{1{\leq}i,j{\leq}r}
=s_{\lambda}(x/y),
\]
while the definition of $\operatorname{sch}$ yields  
\[
\operatorname{sch}(\chi^{\lambda})=\sum_{\mu{\vdash}r}\chi^{\lambda}(\mu)
\dfrac{p_{\mu}(x/y)}{\mu?}.
\]
\end{proof}
Let 
\[
v_m(q)=\prod_{i=1}^m\dfrac{1-q^i}{1-q}
\]
and
\[
v_{\lambda}(q)=\prod_{i=1}^lv_{\lambda_i}(q)
\]
for $\lambda=(\lambda_1,\lambda_2,{\hdots}\lambda_m){\vdash}r,\,l(\lambda){\leq}m$, 
(in which some of the $\lambda_i$ may be zero). 
The Hall-Littlewood symmetric function $P_{\lambda}(x;q)$ for the variables 
$x=(x_1,x_2,{\hdots},x_m)$ is defined by (\cite{Mac},III,(2.1))
\[
P_{\lambda}(x;q)=\dfrac{1}{v_{\lambda}(q)}\sum_{w{\in}\mathfrak{S}_m}
w\Big{(}x_1^{\lambda_1}x_2^{\lambda_2}{\cdots}x_m^{\lambda_m}\prod_{i<j}\dfrac{x_i-qx_j}{x_i-x_j}\Big{)}.
\]
$P_{\lambda}(x;q)$ is defined to be zero if $l(\lambda)>m$. 
$P_{\lambda}(x;q),\,l(\lambda){\leq}m$ constitute a basis of 
$(\Lambda_m)_{\mathbb{Q}(q)}=\Lambda_m{\otimes}_{\mathbb{Z}}\mathbb{Q}(q)$. 
Let $\varphi_r(q)=\prod_{i=1}^r(1-q^i)$ and 
\[
b_{\lambda}(q)=\prod_{i{\geq}1}{\varphi}_{m_i({\lambda})}(q),
\]
where $m_i({\lambda})={\#}\{j\,|\,\lambda_j=i\}$. 
Another Hall-Littlewood symmetric function $Q_{\lambda}(x;q)$ is defined by (\cite{Mac},III,(2.11))
\[
Q_{\lambda}(x;q)=b_{\lambda}(q)P_{\lambda}(x;q).
\]
$P_{\lambda}(x;q)$ and $Q_{\lambda}(x;q)$ are homogeneous of degree $|{\lambda}|$. 
Furthermore, two other symmetric functions $q_r(x;q)$ and $q_{\lambda}(x;q)$ for a partition  $\lambda=(\lambda_1,\lambda_2,{\hdots},\lambda_{l(\lambda)})$ are defined by 
\begin{equation*}
\begin{split}
q_0(x;q)&=1=P_{0}(x;q) \\
q_r(x;q)&=(1-q)P_{(r)}(x;q)=(1-q)\sum_{i=1}^mx_i^r\prod_{j{\neq}i}\dfrac{x_i-qx_j}{x_i-x_j}{\quad}
(r{\geq}1),\\
q_{\lambda}(x;q)&=\prod_{i=1}^{l(\lambda)}q_{\lambda_i}(x;q).
\end{split}
\end{equation*}
As in \cite{Mac},III,(4.8), we may define a scalar product on $(\Lambda_m)_{\mathbb{Q}(q)}$ by requiring 
that the bases ($q_{\lambda}(x;q)$) and ($m_{\lambda}(x)$) be dual to each other: 
\[
{\langle}q_{\lambda}(x;q),m_{\mu}(x){\rangle}=\delta_{\lambda,\mu}
\]
Then one can see 
\[
{\langle}P_{\lambda}(x;q),Q_{\lambda}(x;q){\rangle}=\delta_{\lambda,\mu}{\quad}\text{(\cite{Mac},III,(4.9))}.
\]
The generating function $Q(u)$ of $q_r(x;q)$ is 
\begin{equation*}
Q(u)=\sum_{r{\geq}0}q_r(x;q)u^r=\prod_{i=1}^m\dfrac{1-x_iqu}{1-x_iu}{\quad}\text{(\cite{Mac},III,(2.10))}.
\end{equation*}
In the similar manner as (2.1), we define a supersymmetric function $q_r(x/y;q){\in}\Lambda_{m,n}$ 
as follows. 
\begin{equation*}
Q_{x/y}(u)=\sum_{r{\geq}0}q_r(x/y;q)u^r=\prod_{i=1}^m\dfrac{1-x_iqu}{1-x_iu}
\prod_{j=1}^n\dfrac{1-y_ju}{1-y_jqu}.
\end{equation*}
From the definition of $q_r(x/y;q)$, we immediately have  
\begin{equation}
q_r(x/y;q)=\sum_{k=0}^rq_k(x;q)q_{r-k}(qy;q^{-1})=\sum_{k=0}^rq^{r-k}q_k(x;q)q_{r-k}(y;q^{-1}).
\end{equation}
In order to define the Hall-Littlewood supersymmetric functions, we shall give a definition 
of the skew Hall-Littlewood symmetric functions $P_{{\lambda}/{\mu}}(x;q)$ according to \cite{Mac},III,5. 
$P_{{\lambda}/{\mu}}(x;q)$ is a symmetric function which is uniquely determined by 
\begin{equation}
{\langle}P_{{\lambda}/{\mu}}(x;q),Q_{\nu}(x;q){\rangle}
={\langle}P_{\lambda}(x;q),Q_{\mu}(x;q)Q_{\nu}(x;q){\rangle}
\end{equation}
More constructive way of definition is as follows. 
For two partitions $\lambda,\mu$ such that $\lambda-\mu=\theta$ is a horizontal strip, 
let $J_{\theta}=\{j{\in}\mathbb{N}\,|\,\theta'_j<\theta'_{j+1}\}$ and 
\[
\psi_{{\lambda}/{\mu}}(q)=\prod_{j{\in}J_{\theta}}(1-q^{m_j({\mu})})
\]
where $m_j({\mu})={\#}\{i\,|\,\mu_i=j\}$. 
For general $\lambda,\mu$, take a semi-standard tableau 
(it is called {\lq\lq}column-strict tableau{\rq\rq} or simply {\lq\lq}tableau{\rq\rq} in \cite{Mac})
$T=(\mu=\lambda_0{\subset}\lambda_1{\subset}{\cdots}{\subset}\lambda_l=\lambda)$ and set 
\[
\psi_T(q)=\prod_{i=1}^l\psi_{{\lambda_i}/{\lambda_{i-1}}}(q).
\]
Then $P_{{\lambda}/{\mu}}(x;q)$ is defined by 
\[
P_{{\lambda}/{\mu}}(x;q)=\sum_T\psi_T(q)x^T
\]
summed over all semi-standard tableau $T$ of shape $\lambda-\mu$. 
\begin{Lem}
For $0<k{\leq}r$, 
\[
P_{(r)/(k)}(x;q)=\sum_{\nu{\vdash}r-k}(1-q)^{l(\nu)}m_{\nu}(x).
\]
\end{Lem}
\begin{proof}
Let $\lambda=(r)$ and $\mu=(k)$, and consider a semi-standard tableau 
$T=(\mu=\lambda_0{\subset}\lambda_1{\subset}{\cdots}{\subset}\lambda_l=\lambda)$. 
In each $\psi_{{\lambda_i}/{\lambda_{i-1}}}(q)$, we readily see that $J_{\theta}=\{|\lambda_{i-1}|\}$ 
and $m_{|\lambda_{i-1}|}(\lambda_{i-1})=1$ if $\lambda_{i-1}{\subsetneq}\lambda_i$ and 
$J_{\theta}=\phi$ if $\lambda_{i-1}=\lambda_i$. 
Therefore, for each semi-standard tableau $T$, 
$\psi_T(q)=(1-q)^{{\#}\{i\,|\,\lambda_{i-1}{\subsetneq}\lambda_i\}}$. 
Summing up all semi-standard tableaux, we obtain the equation. 
\end{proof}
\begin{Lem}
\begin{equation*}
\begin{split}
&P_{(r)/0}(x;q)=P_{(r)}(x;q), \\
&P_{(r)/(r)}(x;q)=1 \\
&P_{(r)/(k)}(x;q)=(1-q)P_{(r-k)}(x;q){\quad}(0<k<r).
\end{split}
\end{equation*}
\end{Lem}
\begin{proof}
The first equation is direct consequence of (2.9). 
The second is obtained from Lemma 2.3. 
For the third, we use \cite{Mac},III.2,Ex.3,
\begin{equation*}
P_{(r)}(x;q)=\sum_{i=0}^{r-1}(-q)^is_{(r-i,1^i)}(x).
\end{equation*}
Since the transition matrix $M(s,m)$ equals to $(K_{{\lambda}{\mu}})$ consisting of Kostka numbers, 
the numbers of semi-standard tableaux of shape $\lambda$ and weight $\mu$, we have 
\[
P_{(r-k)}(x;q)=\sum_{i=0}^{r-k-1}(-q)^i\sum_{\nu{\vdash}r-k}K_{(r-k-i,1^i)\,{\nu}}m_{\nu}(x).
\]
By \cite{Mac},I.6,Ex.2(b), 
\[
K_{(r-k-i,1^i)\,{\nu}}=\binom{l({\nu})-1}{i}.
\]
So, observing that $l({\nu})-1{\leq}r-k-1$, and that if $l({\nu})-1<i$ then $\binom{l({\nu})-1}{i}=0$, 
we obtain 
\[
P_{(r-k)}(x;q)=\sum_{\nu{\vdash}r-k}\sum_{i=0}^{l({\nu})-1}(-q)^i\binom{l({\nu})-1}{i}m_{\nu}(x).
\]
On the other hand, from Lemma 2.3,  
\begin{equation*}
\begin{split}
P_{(r)/(k)}(x;q)&=\sum_{\nu{\vdash}r-k}(1-q)^{l(\nu)}m_{\nu}(x) \\
&=(1-q)\sum_{\nu{\vdash}r-k}(1-q)^{l({\nu})-1}m_{\nu}(x) \\
&=(1-q)\sum_{\nu{\vdash}r-k}\Big{\{}\sum_{i=0}^{l({\nu})-1}\binom{l({\nu})-1}{i}(-q)^i\Big{\}}m_{\nu}(x).
\end{split}
\end{equation*}
\end{proof}
\begin{Lem}
For $r{\geq}1$, 
\[
{\omega}P_{(r)}(x;q)=(-1)^{r-1}q^{-1}P_{(r)}(qx;q^{-1})
\]
\end{Lem}
\begin{proof}
From (2.25), we have 
\begin{equation*}
\begin{split}
{\omega}P_{(r)}(x;q)&={\omega}\sum_{i=0}^{r-1}(-q)^is_{(r-i,1^i)}(x) \\
&=\sum_{i=0}^{r-1}(-q)^is_{(r-i,1^i)'}(x) \\
&=\sum_{i=0}^{r-1}(-q)^is_{(i+1,1^{r-i-1})}(x) \\
&=\sum_{j=0}^{r-1}(-q)^{r-1-j}s_{(r-j,1^j)}(x) \\
&=(-q)^{r-1}\sum_{j=0}^{r-1}(-q)^{-j}s_{(r-j,1^j)}(x) \\
&=(-q)^{r-1}P_{(r)}(x;q^{-1}) \\
&=(-1)^{r-1}q^{-1}P_{(r)}(qx;q^{-1})
\end{split}
\end{equation*}
\end{proof}
It is known (\cite{Mac},III,(5.5')) that 
\begin{equation*}
P_{\lambda}(x,y;q)=\sum_{\mu{\subset}\lambda}P_{\mu}(x;q)P_{{\lambda}/{\mu}}(y;q).
\end{equation*}
We define the Hall-Littlewood supersymmetric functions $P_{\lambda}(x/y;q)$ to be 
\begin{equation}
P_{\lambda}(x/y;q)
={\omega_y}P_{\lambda}(x,-y;q)
=\sum_{\mu{\subset}\lambda}P_{\mu}(x;q){\omega}P_{{\lambda}/{\mu}}(-y;q).
\end{equation}
When $q=0$, R.H.S. of (2.10) reduces to that of (2.3). Hence we obtain 
$P_{\lambda}(x/y;0)=s_{\lambda}(x/y)$. 
\begin{Prop}
\begin{equation*}
\begin{split}
&P_0(x/y;q)=q_0(x/y;q) \\
&P_{(r)}(x/y;q)=\dfrac{1}{1-q}q_r(x/y;q){\quad}(r>0).
\end{split}
\end{equation*}
\end{Prop}
\begin{proof}
The first equation is obvious. 
From Lemma 2.4, Lemma 2.5 and (2.10), for $r>0$ we have 
\begin{equation*}
\begin{split}
P_{(r)}(x/y;q)&=\sum_{k=0}^rP_{(k)}(x;q){\omega}P_{(r)/(k)}(-y;q) \\
&={\omega}P_{(r)/(0)}(-y;q)+\sum_{k=1}^{r-1}P_{(k)}(x;q){\omega}P_{(r)/(k)}(-y;q)+P_{(r)}(x;q) \\
&=(-1)^r{\omega}P_{(r)}(y;q)+\sum_{k=1}^{r-1}P_{(k)}(x;q)(-1)^{r-k}(1-q){\omega}P_{(r-k)}(y;q)
+P_{(r)}(x;q) \\
&=(-q)^{-1}P_{(r)}(qy;q^{-1})
+\sum_{k=1}^{r-1}P_{(k)}(x;q)(1-q^{-1})P_{(r-k)}(qy;q^{-1})+P_{(r)}(x;q) \\
&=\dfrac{1}{1-q}q_r(qy;q^{-1})
+\sum_{k=1}^{r-1}\dfrac{1}{1-q}q_k(x;q)q_{r-k}(qy;q^{-1})+\dfrac{1}{1-q}q_{r}(x;q) \\
&=\dfrac{1}{1-q}q_r(x/y;q).
\end{split}
\end{equation*}
\end{proof}
\section{Iwahori-Hecke algebras and general quantum super algebras}
In this section, we review the sign $q$-permutation representation of 
the Iwahori-Hecke algebra $\mathscr{H}_q$ of type $A$ (\cite{Mit2}), and the vector representation 
of the general quantum super algebra $U_q^\sigma$ (\cite{B-K-K},\cite{Mit2}). \par
Let $R$ be a commutative domain with $1$, and let $q$ be an invertible 
element of $R$. The Iwahori-Hecke algebra $\mathscr{H}_{R,r}(q)$ of type $A$ is an $R$-algebra 
generated by $\{T_i\,|\,i=1,2,\hdots,r-1\}$ with the relations: 
\renewcommand{\theenumi}{\arabic{enumi}}
\renewcommand{\labelenumi}{(H\theenumi)}
\begin{enumerate}
\item
$T_i^2 = (q-q^{-1})T_i+1$ \qquad if $i=1,2,\hdots,r-1$,
\item
$T_iT_{i+1}T_i = T_{i+1}T_iT_{i+1}$ \qquad if $i=1,2,\hdots,r-2$,
\item
$T_iT_j = T_jT_i$ \qquad if $|i-j|>1$.
\end{enumerate}
\renewcommand{\theenumi}{\arabic{enumi}}
\renewcommand{\labelenumi}{(\theenumi)}
Let $R=\mathbb{Z}[q,q^{-1}]$ be a polynomial ring over $\mathbb{Z}$ with indeterminates $q^{\pm1}$, 
and $K$ the quotient field of $R$. 
Let $M={\oplus}_{k=1}^{m+n}Rv_k$ be a $\mathbb{Z}_2$-graded $R$-module of rank $m+n$. 
By $\mathbb{Z}_2$-graded, we mean that $M$ is a direct sum of two submodule 
$M_{\overline{0}}=\oplus_{k=1}^{m}Rv_k$ and $M_{\overline{1}}=\oplus_{k=m+1}^{m+n}Rv_k$, and that 
for each homogeneous element the degree map $|{\cdot}|$ 
\begin{equation*}
|v|=
\begin{cases}
0& \text{if $v{\in}M_{\overline{0}}$},\\
1& \text{if $v{\in}M_{\overline{1}}$,}
\end{cases}
\end{equation*}
is given. \par
Let $\pi_r$ be the $q$-permutation representation of 
$\mathscr{H}_{R,r}(q)$ 
on the tensor space 
$M^{{\otimes}r}$. As defined in \cite{Mit2}, $\pi_r$ is given by 
$\pi_r(T_i)=\operatorname{Id}^{{\otimes}i-1}{\otimes}T{\otimes}\operatorname{Id}^{{\otimes}r-i-1}$
($i=1,2,\hdots,r-1$) where $T$ is the operator on $M{\otimes}M$ defined by 
\begin{equation*}
Tv_k{\otimes}v_l=
\begin{cases}
\dfrac{(-1)^{|v_k|}(q+q^{-1})+q-q^{-1}}{2}v_k{\otimes}v_l& \text{if $k=l$,}\\
(-1)^{|v_k||v_l|}v_l{\otimes}v_k+(q-q^{-1})v_k{\otimes}v_l& \text{if $k<l$,}\\
(-1)^{|v_k||v_l|}v_l{\otimes}v_k& \text{if $k>l$.}
\end{cases}
\end{equation*}
and $\operatorname{Id}$ is the identity operator on $M$. 
This representation $\pi_r$ is reduced to the (normal) $q$-permutation representation of 
$\mathscr{H}_{R,r}(q)$ with $n=0$ and to the sign permutation representation 
of the symmetric group $\mathfrak{S}_r$ with $q{\rightarrow}1$. 
For abbreviation, we denote $\mathscr{H}_q=\mathscr{H}_{K,r}(q)=K{\otimes_R}\mathscr{H}_{R,r}(q)$ and set 
$V=K{\otimes}_RM$. The above action also defines a representation of $\mathscr{H}_q$ on $V^{{\otimes}r}$. \par
Next, we define the quantum superalgebra $U_q^\sigma\big{(}{\gl}(m,n)\big{)}$ 
and its vector representation on the tensor space $V^{{\otimes}r}$ according to \cite{B-K-K}. 
$U_q^\sigma\big{(}{\gl}(m,n)\big{)}$ is a Hopf algebra obtained from the quantum 
superalgebra $U_q\big{(}{\gl}(m,n)\big{)}$, which is a Hopf superalgebra, 
by adding an involutive element $\sigma$. 
\par
Let $P=\oplus_{b{\in}B}\mathbb{Z}\epsilon_b$ be a free $\mathbb{Z}$-module where 
$B=B_+{\cup}B_-$ with $B_+=\{1,\hdots,m\}$ and $B_-=\{m+1,\hdots,m+n\}$, and 
$\Pi=\{\alpha_i=\epsilon_i-\epsilon_{i+1}\}_{i{\in}I}$ a set of simple roots with the index set 
$I=I_{\rm even}{\cup}I_{\rm odd}$ where 
$I_{\rm even}=\{1,2,\hdots,m-1,m+1,\hdots,m+n-1\}$ and $I_{\rm odd}=\{m\}$. 
We define a map $p:I{\longrightarrow}\{0,1\}$ to be such that
\begin{equation*}
p(i)=
\begin{cases}
0& \text{if $i{\in}I_{\rm even}$},\\
1& \text{if $i{\in}I_{\rm odd}$.}
\end{cases}
\end{equation*}
A $\mathbb{Q}$-valued symmetric bilinear form on $P$ 
$(\cdot,\cdot) : P{\times}P{\longrightarrow}{\mathbb{Q}}$ is defined as follows. 
\[
(\epsilon_a,\epsilon_{a'})=
\begin{cases}
1& \text{if $a=a'{\in}B_+$,}\\
-1& \text{if $a=a'{\in}B_-$,}\\
0& \text{otherwise.}
\end{cases}
\]
The set $\Pi^{\vee}=\{h_i|i{\in}I\}$ of simple coroots is uniquely determined by the formula 
$\ell_i{\langle}h_i,\lambda{\rangle}=(\alpha_i,\lambda)$ for any $\lambda{\in}P$, where 
${\langle}\cdot,\cdot{\rangle}$ is the natural pairing 
${\langle}\cdot,\cdot{\rangle} : P^*{\times}P{\longrightarrow}{\mathbb{Z}}$ between $P$ and $P^*$ and 
\[
\ell_i=
\begin{cases}
1& \text{if $i=1,\hdots,m$,}\\
-1& \text{if $i=m+1,\hdots,m+n-1$.}
\end{cases}
\]
The quantized enveloping algebra $U_q^\sigma\big{(}{\gl}(m,n)\big{)}$ is the unital associative algebra over 
$K$ with generators $q^h (h{\in}P^*),e_i,f_i (i{\in}I)$ and an additional element 
$\sigma$ which satisfy the following defining relations:
\renewcommand{\theenumi}{\arabic{enumi}}
\renewcommand{\labelenumi}{(Q\theenumi)}
\begin{enumerate}
\item
$q^h=1$ \quad for $h=0$, 
\item
$q^{h_1}q^{h_2}=q^{h_1+h_2}$ \quad for $h_1,h_2{\in}P^*$, 
\item
$q^he_i=q^{{\langle}h,\alpha_j{\rangle}}e_iq^h$ \quad for $h{\in}P^*$ and $i{\in}I$,
\item
$q^hf_i=q^{-{\langle}h,\alpha_j{\rangle}}f_iq^h$ \quad for $h{\in}P^*$ and $i{\in}I$,
\item
$[e_i,f_j]=\delta_{ij}\dfrac{q^{\ell_ih_i}-q^{-\ell_ih_i}}{q^{\ell_i}-q^{-\ell_i}}$ 
\quad for $i,j{\in}I$,
\item
$\sigma^2=1$,
\item
$q^h{\sigma}={\sigma}q^h$ \quad for $h{\in}P^*$,
\item
$e_i{\sigma}=(-1)^{p(i)}{\sigma}e_i$ \quad for $i{\in}I$,
\item
$f_i{\sigma}=(-1)^{p(i)}{\sigma}f_i$ \quad for $i{\in}I$, 
\end{enumerate}
\renewcommand{\theenumi}{\arabic{enumi}}
\renewcommand{\labelenumi}{(\theenumi)}
where $[e_i,f_j]$ means the supercommutator 
\begin{equation*}
[e_i,f_j]=e_if_j-(-1)^{p(i)p(j)}f_je_i.
\end{equation*}
We assume further conditions (bitransitivity condition, see \cite{Kac} p.19):
\renewcommand{\theenumi}{\arabic{enumi}}
\renewcommand{\labelenumi}{(Q\theenumi)}
\begin{enumerate}
\item[(Q10)]
If $a{\in}\sum_{i{\in}I}U_q(\mathfrak{n}_+)e_iU_q(\mathfrak{n}_+)$ satisfies 
$f_ia{\in}U_q(\mathfrak{n}_+)f_i$ for all $i{\in}I$, then $a=0$, 
\item[(Q11)]
If $a{\in}\sum_{i{\in}I}U_q(\mathfrak{n}_-)f_iU_q(\mathfrak{n}_-)$ satisfies 
$e_ia{\in}U_q(\mathfrak{n}_-)e_i$ for all $i{\in}I$, then $a=0$, 
\end{enumerate}
\renewcommand{\theenumi}{\arabic{enumi}}
\renewcommand{\labelenumi}{(\theenumi)}
where $U_q(\mathfrak{n}_+)$ (respectively $U_q(\mathfrak{n}_-)$) is the subalgebra of 
$U_q^\sigma\big{(}{\gl}(m,n)\big{)}$ generated by $\{e_i|i{\in}I\}$ (respectively $\{f_i|i{\in}I\}$). 
$U_q^\sigma\big{(}{\gl}(m,n)\big{)}$ is a Hopf algebra whose comultiplication $\triangle$, 
counit $\varepsilon$, antipode $S$ are as follows. 
\begin{equation*}
\begin{split}
&\triangle(\sigma)=\sigma{\otimes}\sigma, \\
&\triangle(q^h)=q^h{\otimes}q^h \quad \text{for $h{\in}P^*$}, \\
&\triangle(e_i)=e_i{\otimes}q^{-\ell_ih_i}+\sigma^{p(i)}{\otimes}e_i \quad 
\text{for $i{\in}I$}, \\
&\triangle(f_i)=f_i{\otimes}1+\sigma^{p(i)}q^{\ell_ih_i}{\otimes}f_i \quad 
\text{for $i{\in}I$}, \\
&\varepsilon(\sigma)=\varepsilon(q^h)=1 
\quad \text{for $h{\in}P^*$}, \quad 
\varepsilon(e_i)=\varepsilon(f_i)=0 \quad \text{for $i{\in}I$}, \\
&S(\sigma)=\sigma, \quad S(q^{{\pm}h})=q^{{\mp}h} 
\quad \text{for $h{\in}P^*$}, \\
&S(e_i)=-\sigma^{p(i)}e_iq^{\ell_ih_i},{\quad}S(f_i)
=-\sigma^{p(i)}q^{-\ell_ih_i}f_i \quad \text{for $i{\in}I$}. 
\end{split}
\end{equation*}
For the sake of abbreviation, we denote $U^\sigma_q\big{(}\gl(m,n)\big{)}$ by $U^\sigma_q$. 
The vector representation ($\rho,V$) of $U^\sigma_q$ on $\mathbb{Z}_2$-graded 
vector space $V=V_{\bar{0}}{\oplus}V_{\bar{1}}$ (recall that 
$V_{\bar{0}}=\oplus_{i=1}^{m}Rv_i,V_{\bar{1}}=\oplus_{i=m+1}^{m+n}Rv_i$) is defined by 
\begin{equation*}
\begin{split}
&\rho(\sigma)v_j=(-1)^{|v_j|}v_j \quad \text{for $j=1,\hdots,m+n$}, \\
&\rho(q^h)v_j=q^{{\langle}h,\epsilon_j{\rangle}}v_j \quad \text{for $h{\in}P^*,j=1,\hdots,m+n$}, \\
&\rho(e_j)v_{j+1}=v_j \quad \text{for $j=1,\hdots,m+n-1$}, \\
&\rho(f_j)v_j=v_{j+1} \quad \text{for $j=1,\hdots,m+n-1$}, \\
&\text{otherwise $0$.}
\end{split}
\end{equation*}
This representation can be extended to the representation on the tensor space $V^{{\otimes}r}$. 
Let ${\rho}_r$ be the map from $U^\sigma_q$ to 
$\operatorname{End}_{K}(V^{{\otimes}r})$ defined by 
\begin{equation*}
\begin{split}
&{\rho}_r(\sigma)=\rho(\sigma)^{{\otimes}r}, \\
&{\rho}_r(q^h)=\rho(q^h)^{{\otimes}r} \quad \text{for $h{\in}P^*$}, \\
&{\rho}_r(e_i)=\sum_{k=1}^N\rho(\sigma^{p(i)})^{{\otimes}k-1}{\otimes}\rho(e_i) 
{\otimes}\rho(q^{-\ell_ih_i})^{{\otimes}r-k} \quad 
\text{for $i{\in}I$}, \\
&{\rho}_r(f_i)=\sum_{k=1}^r\rho(\sigma^{p(i)}q^{\ell_ih_i})^{{\otimes}k-1}{\otimes}\rho(f_i)
{\otimes}\operatorname{Id}^{{\otimes}r-k} \quad 
\text{for $i{\in}I$}.
\end{split}
\end{equation*}
Let $\triangle^{(1)}=\triangle$ and set 
$\triangle^{(k)}=(\triangle{\otimes}\operatorname{Id}^{{\otimes}k-1})\triangle^{(k-1)}$ 
inductively. Then from the definition of $\triangle$, we have 
${\rho}_r(g)=\rho^{{\otimes}r}\circ\triangle^{(r-1)}(g)$ for $g{\in}U^\sigma_q$ 
immediately. 
\begin{Prop}[\cite{B-K-K} Proposition 3.1]
${\rho}_r$ gives a completely reducible representation of $U^\sigma_q$ 
on $V^{{\otimes}r}$ for $r{\ge}1$. 
\end{Prop}
We denote $\pi_r\big{(}\mathscr{H}_q\big{)}$ by $\mathcal{A}_q$ and 
$\rho_r\big{(}U^\sigma_q\big{)}$ by $\mathcal{B}_q$. 
In \cite{Mit2}, we have shown that $\mathcal{A}_q$ and $\mathcal{B}_q$ are full centralizers 
of each other in $\operatorname{End}_{K}V^{{\otimes}r}$. Namely, 
\begin{Theo}[\cite{Mit2} Theorem 4.4]
$\operatorname{End}_{\mathcal{B}_q}V^{{\otimes}r}
=\mathcal{A}_q$ and $\operatorname{End}_{\mathcal{A}_q}V^{{\otimes}r}
=\mathcal{B}_q$. 
\end{Theo}
We denote by $\overline{K}$ the algebraic closure of $K$.
We set 
$\overline{\mathscr{H}}_q=\mathscr{H}_q{\otimes_{K}}\overline{K}$,
$\overline{U}^\sigma_q=U^\sigma_q{\otimes_{K}}\overline{K}$, 
$\overline{\mathcal{A}}_q=\mathcal{A}_q{\otimes_{K}}\overline{K}$, 
$\overline{\mathcal{B}}_q=\mathcal{B}_q{\otimes_{K}}\overline{K}$. 
Then, 
$\pi_r\big{(}\overline{\mathscr{H}}_q\big{)}=\overline{\mathcal{A}}_q$ and 
$\rho_r\big{(}\overline{U}^\sigma_q\big{)}=\overline{\mathcal{B}}_q$ 
as $\overline{K}$-algebras of operators on 
$\overline{V}^{{\otimes}r}=(V{\otimes_{K}}\overline{K})^{{\otimes}r}$. 
Let $H(m,n;r)=\{\lambda=(\lambda_1,\lambda_2,\hdots){\vdash}r|\lambda_j{\leq}n 
\text{ if } j>m\}$. Diagrams of elements of $H(m,n;r)$ are exactly those contained 
by the $(m,n)$-hooks. 
\begin{Theo}[\cite{Mit2} Theorem 5.1]
$\overline{V}^{{\otimes}r}=\bigoplus_{\lambda{\in}H(m,n;r)}H_{\lambda}{\otimes}V_{\lambda}$ where, \\
$H_{\lambda}\,(\lambda{\in}H(m,n;r))$ are mutually non-isomorphic simple 
$\overline{\mathscr{H}}_q$-modules, \\
$V_{\lambda}\,(\lambda{\in}H(m,n;r))$ are mutually non-isomorphic simple 
$\overline{U}^\sigma_q$-modules.
\end{Theo}
In fact, Theorem 3.3 holds over $K$ since it is known that $K$ is a splitting field for $\mathscr{H}_q$. 
Hence $\mathscr{H}_q$ is a split semisimple $K$-algebra, and we do not need to take the algebraic closure. 
\section{Traces of the actions of $\mathscr{H}_q$ and $U^\sigma_q$}
In this section, we introduce indeterminates $x_1,{\hdots},x_m,y_1,{\hdots},y_n$ associated to 
the basis $v_1,{\hdots},v_{m+n}$ and consider the trace of the product of two operators 
$D_r{\in}\mathcal{B}_q$ and $\pi_r(h){\in}\mathcal{A}_q$ with 
$h{\in}\mathscr{H}_q$, 
in the same manner as in \cite{Ram} or \cite{Shoji}. Using a partition of unity of 
$\mathscr{H}_q$ which has been given in \cite{Gyoja} and \cite{Wenzl} independently, 
we obtain that the trace generates characters of $\mathscr{H}_q$. 
Regarding $\mathbb{C}[\mathfrak{S}_r]$ as the specialization of $\mathscr{H}_q$ to $1$, 
we use same notations such as $\pi_r,\chi^{\lambda}$,etc., for $\mathbb{C}[\mathfrak{S}_r]$. \par
Let $K'=K(z_1,z_2,{\hdots},z_{m+n})$ be the field of rational functions on $K$. 
In the remainder of this paper, we assume that $\mathscr{H}_q,U^\sigma_q,V$, etc., are 
defined over $K'$ and use the same notations such as 
$\mathscr{H}_q,U^\sigma_q,V$, etc.. 
We notice that Proposition 3.1, Theorem 3,2, Theorem 3.3 are hold for $K'$. 
Let 
\begin{equation*}
\begin{split}
\mathscr{I}_{r,m+n}&=\big{\{} \mathbf{i}=(i_1,i_2,{\hdots},i_r)\,|\,1{\leq}i_k{\leq}m+n \big{\}} \\
\mathscr{I}_{r,m+n}^+&=\big{\{} \mathbf{i}=(i_1,i_2,{\hdots},i_r)\,|\,1{\leq}i_k{\leq}m+n,\,
i_1{\leq}i_2{\leq}{\cdots}{\leq}i_r \big{\}}
\end{split}
\end{equation*}
and
\[
\mathscr{C}_{r,m+n}=\big{\{} \mathbf{c}=(c_1,c_2,{\hdots},c_{m+n})\,|\,c_k{\geq}0,{\sum}c_k=r \big{\}}.
\]
For $\mathbf{i}=(i_1,i_2,{\hdots},i_r){\in}\mathscr{I}_{r,m+n}$, we define 
$c(\mathbf{i})=(c_1,c_2,{\hdots},c_{m+n})$ where $k$ appears $c_k$ times in $i_1,i_2,{\hdots},i_r$ for 
each $k=1,2,{\hdots},m+n$. 
Clearly $c$ maps $\mathscr{I}_{r,m+n}$ onto $\mathscr{C}_{r,m+n}$. 
Let $E_i{\in}\operatorname{End}_{K'}V\,(i=1,2,{\hdots}m+n)$ be the projections, 
\[
E_iv_k=\delta_{ik}v_k.
\]
We set $v_{\mathbf{i}}=v_{i_1}{\otimes}v_{i_2}{\otimes}{\cdots}{\otimes}v_{i_r}$ and 
\[
E_{\mathbf{c}}=\sum_{c(\mathbf{i})=\mathbf{c}}E_{i_1}{\otimes}{\cdots}{\otimes}E_{i_r}, 
\]
for $\mathbf{c}{\in}\mathscr{C}_{r,m+n}$. 
Then $E_{\mathbf{c}}$ is the projection from $V^{{\otimes}r}$ to the subspace 
$V^{{\otimes}r}_{\mathbf{c}}$ which is defined by 
\[
V^{{\otimes}r}_{\mathbf{c}}=\sum_{c(\mathbf{i})=\mathbf{c}}K'v_{\mathbf{i}}.
\]
By an easy calculation, we see that $E_i{\otimes}E_j+E_j{\otimes}E_i\,(i,j=1,2,{\hdots},m+n)$ 
commute with $T$ on $V^{{\otimes}2}$. 
For each term 
\[
E_{i_1}{\otimes}{\cdots}{\otimes}E_{i_{k-1}}{\otimes}E_{i_k}{\otimes}E_{i_{k+1}}{\otimes}E_{i_{k+2}}
{\otimes}{\cdots}{\otimes}E_{i_r}
\] 
of $E_{\mathbf{c}}$, there exists a term such that 
\[
E_{i_1}{\otimes}{\cdots}{\otimes}E_{i_{k-1}}{\otimes}E_{i_{k+1}}{\otimes}E_{i_k}{\otimes}E_{i_{k+2}}
{\otimes}{\cdots}{\otimes}E_{i_r}
\] 
in $E_{\mathbf{c}}$. Hence $E_{\mathbf{c}}$ commutes with the action of $\mathscr{H}_q$, namely, 
$E_{\mathbf{c}}$ belongs to $\operatorname{End}_{\mathcal{A}_q}V^{{\otimes}r}=\mathcal{B}_q=\rho_r(U^\sigma_q)$. 
Let $z^{\mathbf{c}}=z_1^{c_1}z_2^{c_2}{\cdots}z_{m+n}^{c_{m+n}}$ and define an operator $D_r$ on 
$V^{{\otimes}r}$ by 
\[
D_r=\sum_{\mathbf{c}{\in}\mathscr{C}_{r,m+n}}z^{\mathbf{c}}E_{\mathbf{c}}.
\]
As mentioned above, we have $D_r{\in}\operatorname{End}_{\mathcal{A}_q}V^{{\otimes}r}$. 
In the same manner as in the proof of Lemma 3.5 in \cite{Ram}, one can show: 
\begin{Lem}[\cite{Ram} Lemma 3.5]
For any idempotent $p{\in}\mathscr{H}_q$, $\operatorname{tr}\big{(}D_r\pi_r(p)\big{)}$ is independent of $q$. 
\end{Lem}
Let us define the specialization to a nonzero complex number $t$ to be a ring homomorphism 
$\varphi_t : R{\longrightarrow}\mathbb{C}$ with the condition $\varphi_t(q)=t$. 
$\mathbb{C}$ becomes $(\mathbb{C},R)$-bimodule, with $R$ acting from the right via 
$\varphi_t$. By the specialization $\varphi_t$, one has 
$\mathbb{C}{\otimes}_R\mathscr{H}_{R,r}(q){\cong}\mathscr{H}_{\mathbb{C},r}(t)$, especially, if $t=1$ then 
$\mathbb{C}{\otimes}_R\mathscr{H}_{R,r}(q){\cong}\mathbb{C}[\mathfrak{S}_r]$. 
We denote by $M_t=\mathbb{C}{\otimes}_RM$ the specialization of $M$ by $\varphi_t$. 
If $t$ is a transcendental number, then $K{\cong}\mathbb{Q}(t)$ as fields via $\mathbb{Q}$-homomorphism 
$q\,{\mapsto}\,t$. Therefore, we have $\overline{K}{\cong}\overline{\mathbb{Q}(t)}{\cong}\mathbb{C}$ 
as fields. 
Let $\mathbb{C}'=\mathbb{C}(z_1,z_2,{\hdots},z_{m+n})$. 
Then $\mathscr{H}_{\mathbb{C}',r}(t)$ is a split semisimple $\mathbb{C}'$-algebra. 
We may assume that 
$M_1$ is defined over $\mathbb{C}'$. 
We replace $z_1,z_2,{\hdots},z_m$ by $x_1,x_2,{\hdots},x_m$ and $z_{m+1},z_{m+2},{\hdots},z_{m+n}$ by 
$-y_1,-y_2,{\hdots},-y_n$. Then we have the following two lemmas in a similar way as in the proofs of 
Lemma 3.6 and Lemma 3.7 in \cite{Ram}. 
\begin{Lem}
Let $\gamma_r=s_1s_2{\cdots}s_{r-2}s_{r-1}{\in}\mathfrak{S}_r$. Then 
\[
\operatorname{tr}\big{(}D_r\pi_r(\gamma_r)\big{)}=\sum_{i=1}^mx_i^r+(-1)^{r-1}\sum_{i=1}^{n}(-y_i)^r
=p_r(x)-p_r(y)=p_r(x/y)
\]
where $D_r\pi_r(\gamma_r)$ is regarded as an operator on $M_1^{{\otimes}r}$. 
\end{Lem}
\begin{proof}
Since 
\[
D_r\pi_r(\gamma_r)v_{i_1}{\otimes}v_{i_2}{\otimes}{\cdots}{\otimes}v_{i_r}=
D_r(-1)^{|v_{i_r}|(|v_{i_1}|+|v_{i_2}|+{\cdots}+|v_{i_{r-1}}|)}v_{i_r}{\otimes}v_{i_1}{\otimes}
{\cdots}{\otimes}v_{i_{r-1}},
\]
the coefficient of $v_{i_1}{\otimes}v_{i_2}{\otimes}{\cdots}{\otimes}v_{i_r}$ in 
$D_r\pi_r(\gamma_r)v_{i_1}{\otimes}v_{i_2}{\otimes}{\cdots}{\otimes}v_{i_r}$ is as follows. 
\[
\text{coef.}=
\begin{cases}
x_{i_1}^r{\quad}\text{if $i_1=i_2={\cdots}=i_r{\leq}m$} \\
(-1)^{r-1}(-y_{i_1})^r{\quad}\text{if $i_1=i_2={\cdots}=i_r>m$} \\
0{\quad}\text{otherwise.}
\end{cases}
\]
\end{proof}
\begin{Lem}
Let $p_{\lambda}$ be a minimal idempotent of $\mathbb{C}[\mathfrak{S}_r]$ corresponding to a simple left 
$\mathbb{C}[\mathfrak{S}_r]$-module indexed by $\lambda{\vdash}r$. Then 
\[
\operatorname{tr}\big{(}D_r{\pi}_r(p_{\lambda})\big{)}=s_{\lambda}(x/y)
\]
where $s_{\lambda}(x/y)$ is the supersymmetric Schur function. 
\end{Lem}
\begin{proof}
Let $z_{\lambda}$ be the minimal central idempotent of $\mathbb{C}[\mathfrak{S}_r]$ indexed by 
$\lambda$, and $\chi^{\lambda}$ the character which is afforded by $\mathbb{C}[\mathfrak{S}_r]p_{\lambda}$. 
$z_{\lambda}$ is given by 
\[
z_{\lambda}=\dfrac{d_{\lambda}}{r!}\sum_{\sigma{\in}\mathfrak{S}_r}{\chi^{\lambda}}(\sigma)\sigma.
\]
Let $d_{\lambda}$ be the degree of $\mathbb{C}[\mathfrak{S}_r]p_{\lambda}$. 
$\mathbb{C}[\mathfrak{S}_r]$ contains just $d_{\lambda}$ simple left $\mathbb{C}[\mathfrak{S}_r]$-modules 
which are isomorphic to $\mathbb{C}[\mathfrak{S}_r]p_{\lambda}$. 
Thus we have the following. 
\begin{equation*}
\begin{split}
\operatorname{tr}\big{(}D_r{\pi}_r(p_{\lambda})\big{)}
&=\dfrac{1}{d_{\lambda}}\operatorname{tr}\big{(}D_r{\pi}_r(z_{\lambda})\big{)}\\
&=\dfrac{1}{d_{\lambda}}\dfrac{d_{\lambda}}{r!}\sum_{\sigma{\in}\mathfrak{S}_r}{\chi^{\lambda}}(\sigma)
\operatorname{tr}\big{(}D_r{\pi}_r(\sigma)\big{)}.
\end{split}
\end{equation*}
Let $\mu=(\mu_1,\mu_2,{\hdots}\mu_l){\vdash}r$ where $l=l(\mu)$ is the length of $\mu$. We set 
cycle permutations as follows. 
\begin{equation*}
\begin{split}
\gamma_{\mu_1}&=(1\,2\,{\hdots}\,{\mu_1})\\
\gamma_{\mu_2}&=({\mu_1}+1\,{\mu_1}+2\,{\hdots}\,{\mu_1}+{\mu_2})\\
&\vdots \\
\gamma_{\mu_l}&=({\mu_1}+{\cdots}+{\mu_{l-1}}+1\,{\mu_1}+{\cdots}+{\mu_{l-1}}+2\,{\hdots}\,r)\\
\gamma_{\mu}&=\gamma_{\mu_1}\gamma_{\mu_2}{\cdots}\gamma_{\mu_l}
\end{split}
\end{equation*}
Since $\chi^{\lambda}(\sigma)$ and $\operatorname{tr}\big{(}D_r{\pi}_r(\sigma)\big{)}$ depend 
only upon the cycle type $\mu$ of $\sigma$, they are constant on conjugacy classes $C_{\mu}$. 
Hence, 
\begin{equation*}
\begin{split}
\operatorname{tr}\big{(}D_r{\pi}_r(p_{\lambda})\big{)}
&=\dfrac{1}{r!}\sum_{\mu{\vdash}r}{\chi^{\lambda}(\mu)}
\operatorname{tr}\big{(}D_r{\pi}_r({\gamma}_{\mu})\big{)}|C_{\mu}|\\
&=\dfrac{1}{r!}\sum_{\mu{\vdash}r}{\chi^{\lambda}(\mu)}{\prod_{k=1}^{l(\mu)}}
\operatorname{tr}\big{(}D_r{\pi}_r({\gamma}_{\mu_k})\big{)}|C_{\mu}|\\
&=\dfrac{1}{r!}\sum_{\mu{\vdash}r}{\chi^{\lambda}(\mu)}{\prod_{k=1}^{l(\mu)}}p_{\mu_k}(x/y)|C_{\mu}|\\
&=\sum_{\mu{\vdash}r}{\chi^{\lambda}(\mu)}\dfrac{p_{\mu}(x/y)}{{\mu}?}.
\end{split}
\end{equation*}
From Proposition 2.2, we have $\operatorname{tr}\big{(}D_r{\pi}_r(p_{\lambda})\big{)}=s_{\lambda}(x/y)$ 
as desired. 
\end{proof}
In the same way as in the proof of Theorem 3.8 in \cite{Ram}, we have 
\begin{Theo}
For any $h{\in}\mathscr{H}_q$, 
\[
\operatorname{tr}(D_r{\pi}_r(h))=\sum_{\lambda{\vdash}r}\chi^{\lambda}(h)s_{\lambda}(x/y),
\]
where $\chi^{\lambda}$ is the irreducible character of $\mathscr{H}_q$ corresponding to 
$\lambda$. 
\end{Theo}
\begin{proof}
Let $\{p^{\lambda}_i\,|\,\lambda{\vdash}r,\,i=1,2,{\hdots},d_{\lambda}={\dim}V_{\lambda}\}$ be a partition 
of unity in $\mathscr{H}_q$ such that 
when we apply the specialization $\varphi_1$, $p^{\lambda}_i$ are well-defined and yield 
a partition of unity in $\mathbb{C}[\mathfrak{S}_r]$. Then we have 
\begin{equation*}
\begin{split}
\operatorname{tr}(D_r{\pi}_r(h))&=\sum_{\lambda,\mu{\vdash}r}\sum_{i=1}^{d_{\lambda}}\sum_{j=1}^{d_{\mu}}
\operatorname{tr}(D_r{\pi}_r(p^{\lambda}_ihp^{\mu}_j)) \\
&=\sum_{\lambda,\mu{\vdash}r}\sum_{i=1}^{d_{\lambda}}
\operatorname{tr}(D_r{\pi}_r(p^{\lambda}_ihp^{\lambda}_i)) \\
&=\sum_{\lambda{\vdash}r}\sum_{i=1}^{d_{\lambda}}h^{\lambda}_{ii}\operatorname{tr}(D_r{\pi}_r(p^{\lambda}_i)), 
\end{split}
\end{equation*}
where $h^{\lambda}_{ii}$ is the diagonal element of the representation matrix of $h$ in the irreducible 
representation corresponding to $\lambda$ determined by this partition of unity. 
Thus by Lemma 4.3, we obtain 
\begin{equation*}
\begin{split}
\operatorname{tr}(D_r{\pi}_r(h))&=\sum_{\lambda{\vdash}r}\sum_{i=1}^{d_{\lambda}}h^{\lambda}_{ii}
s_{\lambda}(x/y) \\
&=\sum_{\lambda{\vdash}r}\chi^{\lambda}(h)s_{\lambda}(x/y)
\end{split}
\end{equation*}
\end{proof}
\section{The super Frobenius formula for the characters of $\mathscr{H}_q$}
In this section, we compute $\operatorname{tr}(D_r{\pi}_r(h))$ in detail. 
We shall show that when $h$ is a product of elements corresponding to cycle permutations, 
the trace coincides with a Hall-Littlewood supersymmetric function up to constant. 
There is no notion of cycle type for elements of $\mathscr{H}_q$ in general, 
so not all character values are given in this way. 
Nevertheless, by Ram's result, any character of $\mathscr{H}_q$ is determined by its values on 
elements corresponding to cycle permutations. \par
For $\mathbf{i}=(i_1,i_2,{\hdots}i_r){\in}\mathscr{I}_{r,m+n}$, we define cardinalities of subsets of 
$\{i_1,i_2,{\hdots}i_r\}$ as follows. 
\begin{equation*}
\begin{split}
N^0(\mathbf{i})&={\#}\{j|i_j{\leq}m\},{\quad}N^1(\mathbf{i})={\#}\{j|i_j>m\}, \\
E^0(\mathbf{i})&={\#}\{j|i_j=i_{j+1},i_j{\leq}m\},{\quad}E^1(\mathbf{i})={\#}\{j|i_j=i_{j+1},i_j>m\},{\quad}
E(\mathbf{i})=E^0(\mathbf{i})+E^1(\mathbf{i}), \\
L^0(\mathbf{i})&={\#}\{j|i_j<i_{j+1},i_j{\leq}m\},{\quad}L^1(\mathbf{i})={\#}\{j|i_j<i_{j+1},i_j>m\},{\quad}
L(\mathbf{i})=L^0(\mathbf{i})+L^1(\mathbf{i}).
\end{split}
\end{equation*}
Let $T_{\gamma_r}=T_{s_1}T_{s_2}{\cdots}T_{s_{r-1}}$. 
Using above notations we obtain: 
\begin{Prop}
For $1{\leq}k{\leq}r$, the trace of $D_k{\pi_k}(T_{\gamma_k})$ on $V^{{\otimes}k}$ is given by 
\[
\operatorname{tr}\big{(}D_k{\pi_k}(T_{\gamma_k})\big{)}=\sum_{\mathbf{i}{\in}\mathscr{I}_{k,m+n}^+}
(-1)^{E^1(\mathbf{i})}q^{E^0(\mathbf{i})-E^1(\mathbf{i})}(q-q^{-1})^{L(\mathbf{i})}
z_{i_1}z_{i_2}{\cdots}z_{i_k}
\]
\end{Prop}
\begin{proof}
We prove by induction on $k$. 
If $k=1$, then $T_{\gamma_1}=1$, so the statement holds obviously. 
Now assume $k>1$. We consider 
$D_k{\pi_k}(T_{\gamma_k})v_{i_1}{\otimes}v_{i_2}{\otimes}{\cdots}{\otimes}v_{i_k}$ case-by-case 
depending on the relation between $i_{k-1}$ and $i_k$. \\
case 1 : $i_{k-1}>i_k$
\begin{equation*}
\begin{split}
D_k{\pi_k}(T_{\gamma_k})v_{i_1}{\otimes}v_{i_2}{\otimes}{\cdots}{\otimes}v_{i_{k-1}}
{\otimes}v_{i_k}
=(-1)^{|v_{i_{k-1}}||v_{i_k}|}D_k{\pi_k}(T_{\gamma_{k-1}})v_{i_1}{\otimes}v_{i_2}{\otimes}
{\cdots}{\otimes}v_{i_k}{\otimes}v_{i_{k-1}}
\end{split}
\end{equation*}
Since ${\pi_k}(T_{\gamma_{k-1}})$ acts on the first $k-1$ factors and $i_{k-1}{\neq}i_k$, 
the coefficient of $v_{i_1}{\otimes}v_{i_2}{\otimes}{\cdots}{\otimes}v_{i_{k-1}}{\otimes}v_{i_k}$ 
in $D_k{\pi_k}(T_{\gamma_k})v_{i_1}{\otimes}v_{i_2}{\otimes}{\cdots}{\otimes}v_{i_k}$ is zero. \\
case 2 : $i_{k-1}<i_k$ 
\begin{equation*}
\begin{split}
D_k{\pi_k}(T_{\gamma_k})v_{i_1}{\otimes}v_{i_2}{\otimes}{\cdots}{\otimes}v_{i_{k-1}}
{\otimes}v_{i_k}
&=(-1)^{|v_{i_{k-1}}||v_{i_k}|}D_k{\pi_k}(T_{\gamma_{k-1}})v_{i_1}{\otimes}v_{i_2}{\otimes}
{\cdots}{\otimes}v_{i_k}{\otimes}v_{i_{k-1}}\\
&+(q-q^{-1})D_k{\pi_k}(T_{\gamma_{k-1}})v_{i_1}{\otimes}v_{i_2}{\otimes}
{\cdots}{\otimes}v_{i_{k-1}}{\otimes}v_{i_k}
\end{split}
\end{equation*}
For the same reason as case 1, the coefficient of 
$v_{i_1}{\otimes}v_{i_2}{\otimes}{\cdots}{\otimes}v_{i_{k-1}}{\otimes}v_{i_k}$ in the first term 
is zero. While the one in the second term equals to the coefficient of 
$v_{i_1}{\otimes}v_{i_2}{\otimes}{\cdots}{\otimes}v_{i_{k-1}}$ in 
$(q-q^{-1})z_{i_k}D_{k-1}{\pi_{k-1}}(T_{\gamma_{k-1}})$ $v_{i_1}{\otimes}v_{i_2}{\otimes}
{\cdots}{\otimes}v_{i_{k-1}}$. \\
case 3 : $i_{k-1}=i_k$
\begin{equation*}
\begin{split}
D_k{\pi_k}(T_{\gamma_k})v_{i_1}{\otimes}v_{i_2}{\otimes}{\cdots}{\otimes}v_{i_{k-1}}
{\otimes}v_{i_k}
=\dfrac{(-1)^{|v_{i_k}|}(q+q^{-1})+(q-q^{-1})}{2}D_k{\pi_k}(T_{\gamma_{k-1}})v_{i_1}
{\otimes}v_{i_2}{\otimes}{\cdots}{\otimes}v_{i_{k-1}}{\otimes}v_{i_k}
\end{split}
\end{equation*}
Thus the coefficient of $v_{i_1}{\otimes}v_{i_2}{\otimes}{\cdots}{\otimes}v_{i_k}$ 
in $D_k{\pi_k}(T_{\gamma_k})v_{i_1}{\otimes}v_{i_2}{\otimes}{\cdots}{\otimes}v_{i_k}$ equals to 
that of $v_{i_1}{\otimes}v_{i_2}{\otimes}{\cdots}{\otimes}v_{i_{k-1}}$ in 
$qz_{i_k}D_{k-1}{\pi_{k-1}}(T_{\gamma_{k-1}})$ $v_{i_1}{\otimes}v_{i_2}{\otimes}{\cdots}{\otimes}v_{i_{k-1}}$ 
if $i_k{\leq}m$ and to that of $v_{i_1}{\otimes}v_{i_2}{\otimes}{\cdots}{\otimes}v_{i_{k-1}}$ in 
$-q^{-1}z_{i_k}D_{k-1}$ ${\pi_{k-1}}(T_{\gamma_{k-1}})v_{i_1}{\otimes}v_{i_2}{\otimes}{\cdots}{\otimes}v_{i_{k-1}}$ 
if $i_k>m$. \\
By induction, the assertion follows as desired. 
\end{proof}
Replacing $z_1,z_2,{\hdots},z_m$ by $x_1,x_2,{\hdots},x_m$ and $z_{m+1},z_{m+2},{\hdots},z_{m+n}$ by 
$-y_1,-y_2,{\hdots},-y_n$, we obtain
\begin{equation*}
\begin{split}
\operatorname{tr}\big{(}D_k{\pi_k}(T_{\gamma_k})\big{)}&=\sum_{\mathbf{i}{\in}\mathscr{I}_{k,m+n}^+}
(-1)^{E^1(\mathbf{i})}q^{E^0(\mathbf{i})-E^1(\mathbf{i})}(q-q^{-1})^{L(\mathbf{i})}
z_{i_1}z_{i_2}{\cdots}z_{i_k}\\
&=\sum_{\mathbf{i}{\in}\mathscr{I}_{k,m+n}^+}
q^{E^0(\mathbf{i})}(q-q^{-1})^{L^0(\mathbf{i})}z_{i_1}{\cdots}z_{i_{N^0(\mathbf{i})}}
(-q^{-1})^{E^1(\mathbf{i})}(q-q^{-1})^{L^1(\mathbf{i})}
z_{i_{N^0(\mathbf{i})+1}}{\cdots}z_{i_k} \\
&=\sum_{\mu{\vdash}k}(-q^{-1})^{|\mu|-l({\mu})}(q-q^{-1})^{l({\mu})-1}m_{\mu}(-y)\\
&+\sum_{j=1}^{k-1}\Big{\{}\sum_{\lambda{\vdash}j}q^{|\lambda|-l({\lambda})}(q-q^{-1})^{l({\lambda})}m_{\lambda}(x)
\sum_{\mu{\vdash}k-j}(-q^{-1})^{|\mu|-l({\mu})}(q-q^{-1})^{l({\mu})-1}m_{\mu}(-y) \Big{\}}\\
&+\sum_{\lambda{\vdash}k}q^{|\lambda|-l({\lambda})}(q-q^{-1})^{l({\lambda})-1}m_{\lambda}(x)
\end{split}
\end{equation*}
Let 
\begin{equation*}
\Tilde{q}_k(x;q)=
q^k\sum_{\lambda{\vdash}k}\Big{(}\dfrac{q-q^{-1}}{q}\Big{)}^{l({\lambda})}
m_{\lambda}(x){\quad}\text{if $k>0$},
\end{equation*}
and $\Tilde{q}_0(x;q)=1$. 
Since $m_{\lambda}(x_1,{\hdots},x_m)=0$ if $l(\lambda)>m$, 
$\Tilde{q}_k(x_1,{\hdots},x_m;q)=0$ if $m=0$ and $k>0$. 
From definition of $\Tilde{q}_k(x;q)$, one can describe the trace as follows.  
\begin{equation}
\operatorname{tr}\big{(}D_k{\pi_k}(T_{\gamma_k})\big{)}=
\dfrac{1}{q-q^{-1}}\sum_{j=0}^k\Tilde{q}_j(x;q)\Tilde{q}_{k-j}(-y;-q^{-1}).
\end{equation}
Considering the generating function of $\Tilde{q}_r(x;q)$, we obtain the relation between 
$\Tilde{q}_r(x;q)$ and $q_r(x;q)$ as follows. 
\begin{Lem}
\[
\Tilde{q}_r(x;q)=q_r(qx;q^{-2})=q^rq_r(x;q^{-2})
\]
\end{Lem}
\noindent
This equations are found in \cite{Ram},Theorem4.13 and \cite{Shoji},(6.11.3) in 
slightly different forms, but essentially the same. 
\begin{proof}
Consider the generating function of $\Tilde{q}_r(x;q)$. 
\begin{equation*}
\begin{split}
\sum_{r{\geq}0}\Tilde{q}_r(x;q)u^r&=\sum_{r{\geq}0}q^ru^r
\Big{\{}\sum_{\mu{\vdash}r}\Big{(}\frac{q-q^{-1}}{q}\Big{)}^{l({\mu})}m_{\mu}(x)\Big{\}} \\
&=\prod_{i=1}^m\Big{(}\dfrac{q-q^{-1}}{q}{\cdot}
\dfrac{1}{1-qux_i}-\dfrac{q-q^{-1}}{q}+1\Big{)} \\
&=\prod_{i=1}^m\dfrac{1-q^{-1}ux_i}{1-qux_i} \\
&=\prod_{i=1}^m\dfrac{1-q^{-2}qx_iu}{1-qx_iu} \\
&=\sum_{r{\geq}0}q_r(qx;q^{-2})u^r \\
&=\sum_{r{\geq}0}q^rq_r(x;q^{-2})u^r 
\end{split}
\end{equation*}
Thus the equalities hold. 
\end{proof}
\begin{Theo}
\[
\operatorname{tr}\big{(}D_k{\pi_k}(T_{\gamma_k})\big{)}=\dfrac{q^k}{q-q^{-1}}q_k(x/y;q^{-2})
\]
\end{Theo}
\begin{proof}
From (5.1) and Lemma 5.2, we obtain 
\begin{equation*}
\begin{split}
\operatorname{tr}\big{(}D_k{\pi_k}(T_{\gamma_k})\big{)}&=
\dfrac{1}{q-q^{-1}}\sum_{j=0}^k\Tilde{q}_j(x;q)\Tilde{q}_{k-j}(-y;-q^{-1}) \\
&=\dfrac{1}{q-q^{-1}}\sum_{j=0}^kq_j(qx;q^{-2})q_{k-j}(q^{-1}y;q^2). 
\end{split}
\end{equation*}
On the other hand, by (2.8), we have  
\begin{equation*}
\begin{split}
q_k(x/y;q^{-2})&=
\sum_{j=0}^kq_j(x;q^{-2})q_{k-j}(q^{-2}y;q^2) \\
&=\sum_{j=0}^kq_j(x;q^{-2})q^{j-k}q_{k-j}(q^{-1}y;q^2) \\
&=q^{-k}\sum_{j=0}^kq_j(qx;q^{-2})q_{k-j}(q^{-1}y;q^2).
\end{split}
\end{equation*}
This completes the proof. 
\end{proof}
A finite sequence $\alpha=(\alpha_1,\alpha_2,{\hdots},\alpha_l)$ of positive integers which satisfy 
$\alpha_1+\alpha_2+{\cdots}+\alpha_l=r$ is said to be a composition of $r$ and denote by 
$\alpha{\models}r$. We set cycle permutations for $\alpha$ as follows. 
\begin{equation*}
\begin{split}
\gamma_{\alpha_1}&=(1\,2\,{\hdots}\,{\alpha_1})\\
\gamma_{\alpha_2}&=({\alpha_1}+1\,{\alpha_1}+2\,{\hdots}\,{\alpha_1}+{\alpha_2})\\
&\vdots \\
\gamma_{\alpha_l}&=({\alpha_1}+{\cdots}+{\alpha_{l-1}}+1\,{\alpha_1}+{\cdots}+{\alpha_{l-1}}+2\,{\hdots}\,r)\\
\gamma_{\alpha}&=\gamma_{\alpha_1}\gamma_{\alpha_2}{\cdots}\gamma_{\alpha_l}
\end{split}
\end{equation*}
Then one can readily see that $T_{\gamma_{\alpha}}=T_{\gamma_{\alpha_1}}T_{\gamma_{\alpha_2}}{\cdots}
T_{\gamma_{\alpha_l}}$. 
\begin{Prop}
\[
\operatorname{tr}\big{(}D_r{\pi_r}(T_{\gamma_{\alpha}})\big{)}
=\prod_{i=1}^l\operatorname{tr}\big{(}D_{\alpha_i}{\pi_{\alpha_i}}(T_{\gamma_{\alpha_i}})\big{)}
\]
\end{Prop}
\begin{proof}
It is sufficient to prove for $l=2$. Let $T_{\gamma_{\alpha}}=T_{\gamma_{\alpha_1}}T_{\gamma_{\alpha_2}}$. 
$T_{\gamma_{\alpha_1}}$ and $T_{\gamma_{\alpha_2}}$ act only on 
$v_{i_1}{\otimes}v_{i_2}{\otimes}{\cdots}{\otimes}v_{i_{\alpha_1}}$ and 
$v_{i_{\alpha_1+1}}{\otimes}v_{i_{\alpha_1+2}}{\otimes}{\cdots}{\otimes}v_{i_{\alpha_1+\alpha_2}}$ 
respectively. Therefore, 
\begin{equation*}
\begin{split}
D_r{\pi_r}(T_{\gamma_{\alpha_1}}T_{\gamma_{\alpha_2}})v_{i_1}{\otimes}v_{i_2}{\otimes}{\cdots}{\otimes}v_{i_r}
&=\big{(}D_{\alpha_1}{\pi_{\alpha_1}}(T_{\gamma_{\alpha_1}})
v_{i_1}{\otimes}v_{i_2}{\otimes}{\cdots}{\otimes}v_{i_{\alpha_1}}\big{)}\\
&{\otimes}\big{(}D_{\alpha_2}{\pi_{\alpha_2}}(T_{\gamma_{\alpha_2}})
v_{i_{\alpha_1+1}}{\otimes}v_{i_{\alpha_1+2}}{\otimes}{\cdots}{\otimes}v_{i_{\alpha_1+\alpha_2}}\big{)}.
\end{split}
\end{equation*}
Thus we have proved the assertion. 
\end{proof}
\begin{Theo}
For $\mu{\vdash}r$, 
\[
\dfrac{q^{|{\mu}|}}{(q-q^{-1})^{l(\mu)}}q_{\mu}(x/y;q^{-2})
=\sum_{\lambda{\vdash}r}\chi^{\lambda}(T_{\gamma_{\mu}})s_{\lambda}(x/y)
\]
\end{Theo}
\begin{proof}
The assertion follows from Theorem 5.3 and Proposition 5.4. 
\end{proof}
By Theorem 5.1 in \cite{Ram}, with a slight change for our version, 
we obtain that any character of 
$\mathscr{H}_q$ is determined by its values on $T_{\gamma_{\mu}},\,\mu{\vdash}r$. 
\begin{Theo}[\cite{Ram},Theorem 5.1]
For each $T_{\sigma},\,\sigma{\in}\mathfrak{S}_r$, there exists a $\mathbb{Z}[q,q^{-1}]$ linear combination 
\[
c_{\sigma}=\sum_{\mu{\vdash}r}a_{{\sigma}{\mu}}T_{\gamma_{\mu}},
\]
$a_{{\sigma}{\mu}}{\in}\mathbb{Z}[q,q^{-1}]$, such that $\chi(T_{\sigma})=\chi(c_{\sigma})$ 
for all characters $\chi$ of $\mathscr{H}_q$. 
\end{Theo}

\end{document}